# The independence of the parallel postulate and the acute angles hypothesis on a two right-angled isosceles quadrilateral

Christos Filippidis[1], Prodromos Filippidis

We trace the development of arguments for the consistency of non-Euclidean geometries and for the independence of the parallel postulate, showing how the arguments become more rigorous as a formal conception of geometry is introduced. We use these philosophical views to explain why the certainty of Euclidean geometry was threatened by the development of what we regard as alternatives to it. Finally, we provide a model that creates the basis for proof and demonstration using natural language in order to study the acute angles hypothesis on a two right-angled isosceles quadrilateral, with a rectilinear summit side.

Keywords: parallel postulate; Hilbert axioms; consistency; Saccheri quadrilateral; function theory;

## The independence of the parallel postulate and development of consistency proofs

Euclid's parallel postulate states: "If a straight line falling on two straight lines makes the interior angles on the same side less than two right angles, the two straight lines if produced indefinitely meet on that side on which are the angles less than two right angles." The earliest commentators found fault with this statement as being not self-evident 'Lewis [1920]'. Euclid uses terms such as "indefinitely" and makes logical assumptions that had not been proven or stated. Thus the parallel postulate seemed less obvious than the others.According to Lambert, Euclid overcomes doubt by means of postulates. Euclid's theory thus owes its justification not to the existence of the surfaces that satisfy it, but to the postulates according to which these "models" are constructed ' Dunlop [2009]'.

---

[1] cfjs@outlook.com

For two thousand years, many attempts were made to prove the parallel postulate by using Euclid's first four postulates. Finally the mathematical society have proven the consistency of non-euclidean plane geometry and the independence of the euclidean parallel postulate; and extend that in stating that it would be impossible to prove Euclid's parallel postulate from the other assumptions made by him, since this would involve the denial of the parallel postulate of Bolyai and Lobatschewsky 'CARSLAW [1910]'.

These conclusions mainly derived from Beltrami's (1868) seminal papers concerning the existence of non-Euclidean objects and by the arguments offered by Jules Houel in 1860-1870 for the unprovability of the parallel postulate and for the existence of non-Euclidean geometries. The final step towards rigorous consistency proofs is taken in the 1880s by Henri Poincare.

What Houel sees in Beltrami's putative demonstration that regions of pseudo-spheres are complete Bolyai-Lobatschewsky planes is the presentation of an instance in which all of the properties of the euclidean plane are present except for the property of unique parallelism 'Scanlan [1988]'. A trivial example will make the principle evident.

"If we wish to show that the property of upright posture in mammals does not imply intelligence (in this sense of 'implication'), it is sufficient to show that there are some mammals in which upright posture is present without intelligence (kangaroos), even though there may be others in which upright posture occurs together with intelligence (humans). The application of the result here to the human case, e.g. if we want to say that upright posture in humans does not imply their intelligence, turns on referring to the same properties in the case of both humans and kangaroos" 'Scanlan [1988]'.

It's obvious that the independence of the euclidean parallel postulate is well formed, but a detailed looked at the writings will show that some claims may are overstated or misinterpreted. The most striking example of a non-logical use of the word "independent" is in Bolyai's title of his famous appendix, one of the founding texts of hyperbolic geometry, Absolute Geometry: Independent of the Truth or Falsity of Euclids Axiom XI (which can never be decided a priori).

Bolyai does not, of course, claim that the parallel postulate is independent from the rest of Euclidean geometry in the sense required to show that the parallel postulate is unprovable (though he does claim that it cannot be proven). He means simply that he can prove some theorems of geometry without relying on either the parallel postulate or its negation 'Stump [2007]'.

Sometimes, it was simply stated that the attempt to prove the parallel postulate had gone on too long and that repeated lack of success shows that it is impossible to prove the parallel postulate. On the other hand the geometers working in order to support the existence of non-Euclidean objects were in the pursue of some very impressive results. That's why Beltrami's work had a profound impact on Houel. Houel immediately after he saw it, announced that Beltrami had shown that it is impossible to prove the Parallel Postulate. Strikingly, he made this announcement in eight different journals. In reality they wanted to strongly support and consolidate the initial hypothesis of the acute angles.

The geometers had made just another logical assumption, believing that by finding a contradiction to Saccheri's quadrilateral would weaken the foundations of non-Euclidean geometries. We can also realize that some statements wasn't very well formed or may were just some logical assumptions and we have to treat them with the same way we did with the logical assumptions that Euclid did. 'H. S. CARSLAW [1910]' states: "How ever far the hyperbolic geometry were developed, no contradictory results could be obtained". Stating "How ever far" is equivalent with the term "indefinitely" that Euclid uses.

**The acute angles hypothesis on a two right-angled isosceles quadrilateral**

We consider a Hilbert plane, satisfying the axioms (I1)-(I3), (B1)-(B4), (C1)-(C6), and the axiom of continuity, with the definitions and theorems deduced by them 'Hartshorne'. Based on these considerations we provide a model that creates the basis for proof and demonstration using natural language. Within this model we study the acute angles hypothesis on a two right-angled isosceles quadrilateral, with a rectilinear summit side.

(1) We construct a Saccheri rectilinear quadrilateral, the hypothesis of the acute angles is true.

(2) We form a rectilinear isosceles triangle by using a point on the line of the common perpendicular of the Saccheri rectilinear quadrilateral (the common perpendicular is defined by the base and the summit side of the quadrilateral), and the summit side of the Saccheri rectilinear quadrilateral.

(3) By construction, the perpendicular bisectors to the congruent sides of the isosceles triangle as well as the bisector of the vertex angle of the isosceles triangle are perpendiculars to the base of the acute-angled Saccheri rectilinear quadrilateral. Since the perpendicular to a line through a point not on the line is unique, we can conclude that they are not intersecting.

(4) We consider the cartesian system where the x-axis coincides with the base of the above isosceles triangle and the y-axis coincides with the bisector of the vertex angle of the triangle.

(5) The function which represents the straight line of each perpendicular bisector to the congruent sides of the isosceles triangle is defined by the two points of the perpendicular bisector intersecting the base and the corresponding congruent side of this isosceles triangle.

(6) These two functions representing the two perpendicular bisectors to the congruent sides of the isosceles triangle could be either equal, foreign, or intersecting.

(7) These two functions which represent the two perpendicular bisectors to the congruent sides of the isosceles triangle assign not-equal arguments (x-values) to equal values (y-value). Since the line formed by the two points of intersection to each perpendicular bisector with the corresponding congruent side of the isosceles triangle is perpendicular to y-axis 'Filippidis, P [2016]'.

(8) Since the functions which represent the two perpendicular bisectors to the congruent sides of the isosceles, are not equal neither foreign they are intersecting.

(9) The point of intersection of the perpendicular bisectors to the congruent sides of the isosceles triangle lies on the y-axis, which coincides with the bisector of the vertex angle of this triangle.

(10) The hypothesis that the summit angles of a rectilinear Saccheri quadrilateral are acute is not true.

(11) We can follow the same logical steps regarding the obtuse angles hypothesis on a Saccheri rectilinear quadrilateral.

(12) The type of the summit side angles (acute, obtuse, right) of a Saccheri quadrilateral depends on the curvature of the summit side.